\documentclass[11pt,bezier]{article}
\usepackage{amsmath,amssymb,amsfonts,euscript}
\usepackage[usenames]{color}
\newpage

\newcommand{\be}{\begin{equation}}
\newcommand{\ee}{\end{equation}}
\newtheorem{theorem}{Theorem}[section]

\newtheorem{property}{Property}[section]
\newtheorem{lemma}{Lemma}[section]

\newtheorem{example}{Example}[section]

\title{\bf\Large Two New Entropy Estimators for Testing Exponentiality with Type-II Censored Data
}
\vspace{-1cm}
\author
{A.Kohansal., S.Rezakhah \footnote{Faculty of Mathematics and Computer Science, Amirkabir University of Technology,
Tehran, Iran. Email: rezakhah@aut.ac.ir, ak\_kohansal@aut.ac.ir}}
\begin{document}
\maketitle

\begin{abstract}
This paper proposes two estimators of the joint entropy of the Type-II censored data. Consistency of both estimators is proved. Simulation results show that the second one shows less bias and root of mean square error (RMSE) than leading estimator. Also, two goodness of fit test statistics based on the Kullback-Leibler information with the Type-II censored data are established and their performances with the leading test statistics are compared. We provide a Monte Carlo simulation study which shows that the test statistics $T^{(1)}_{m,n,r}$ and $T^{(2)}_{m,n,r}$ show better powers than leading test statistics against the alternatives with monotone decreasing and monotone increasing hazard functions, respectively.\\ \quad \\
{\it Keywords}: Entropy, Monte Carlo simulation, Kullback-Leibler information, Moving average method, Hazard function. \\ \quad \\
{\it Mathematics Subject Classification:} 62G10, 62G30.
\end{abstract}
\section{Introduction}
\par
Suppose that a random variable $X$ has a distribution function $F(x)$, with a continuous density function $f(x)$. The differential entropy $H(f)$ of the random variable is defined by Shannon \cite{sha} to be
\begin{equation}\label{ent}
H(f)=-\int_{-\infty}^{\infty}f(x)\ln f(x)\mathrm{d}x.
\end{equation}
The entropy difference $H(f)-H(g)$ has been considered in \cite{dud2} and \cite{gok} for establishing the goodness of fit tests for the class of the maximum entropy distributions.
\par
The Kullback-Leibler (KL) information in favor of $g(x)$ against $f(x)$ is defined to be
$$I(g;f)=\int_{-\infty}^{\infty}g(x)\ln\frac{g(x)}{f(x)}\mathrm{d}x.$$
Because $I(g;f)$ has the property that $I(g;f)\geq0$, and the equality holds if and only if $f=g$, different estimators of the KL information has been also considered as a goodness of fit test statistic in some papers including \cite{ari}, \cite{ebr2}, \cite{par3} and \cite{zam1}. For complete samples, some of these test statistics perform very well for exponentiality \cite{ebr2}, and some others of them perform very well for normality, see \cite{vas}, \cite{zam2} and \cite{nou}.
\par
For the censored data, some authors studied the problem of goodness of fit test and discussed some test statistics. Brain and Shapiro \cite{bra} proposed two test statistics and show that these test statistics perform better than other test statistics for the censored data. Samanta and Schwarz \cite{sam} proposed a test statistic and showed that the proposed test statistic has competing performance with the test statistics which introduced by Brain and Shapiro \cite{bra} for the censored data. Recently Park \cite{par2} obtained an estimator for entropy of Type-II censored data and proposed a test statistic based on KL information. He showed that the power of the proposed test statistic is greater than the power of the test statistics which proposed by Brain and Shapiro \cite{bra}, and Samanta and Schwarz \cite{sam} against the alternatives with monotone increasing hazard functions. In the case of progressively censored data, Balakrishnan et al. \cite{bala} studied the testing exponentiality based on KL information with progressively Type-II censored data. Habibi Rad et al. \cite{hab} studied goodness of fit test based on KL information for progressively Type-II censored data. Pakyari and Balakrishnan \cite{pak1} proposed several goodness of fit methods for location-scale families of distributions under progressively Type-II censored data. They \cite{pak2} also investigated a general purpose approximate goodness of fit test for progressively Type-II censored data.
\par
In this paper, we enhance the estimator which was introduced by Park \cite{par2} and obtain two new entropy estimators of Type-II censored data. Simulation results show that the second one shows less bias and RMSE than leading estimator. Also, we provide two new test statistics. The first one achieves higher power than the previous test statistics against the alternatives with monotone decreasing hazard functions and the other one achieves higher power than the previous test statistics against the alternatives with monotone increasing hazard functions.
\par
The rest of the article is arranged as follows: In Section 2, we introduce two estimators of the joint entropy of the Type-II censored data. Also, we show that both are consistent. Scale invariance property of variances and mean squared errors of the proposed estimators is studied in the same section. In Section 3, we use the KL information with the Type-II censored data and obtain two new test statistics. In Section 4, we introduce goodness-of-fit tests for exponentiality based on the proposed test statistics and then compare their powers with the powers of other test statistics. Also, by using the new test statistics, we compare biases and RMSEs of the new entropy estimators with the leading entropy estimator.

\section{New entropy estimators}
In this section, we introduce two entropy estimators and prove some of their properties.

\subsection{Entropy estimator for monotone decreasing hazard function alternatives}
\par
In this subsection, we obtain one entropy estimator which provides a new test statistic that achieves higher power than the previous test statistics against the alternatives with monotone decreasing hazard functions. \par
Vasicek \cite{vas} expressed (\ref{ent}) in the form,
\begin{equation*}
H=\int_{0}^{1}\ln\left(\frac{\mathrm{d}F^{-1}p}{\mathrm{d}p}\mathrm{d}p\right)
\end{equation*}
and provided its estimator as: $$\displaystyle H(m,n)=\frac{1}{n}\sum_{i=1}^n\ln\frac{x_{(i+m:n)}-x_{(i-m:n)}}{\frac{2m}{n}},$$ where
the window size $m$ is a positive integer, which is less than $n/2$; and $x_{(i:n)}=x_{(1:n)}$ for $i<1$, and $x_{(i:n)}=x_{(n:n)}$ for $i>n$. Recently Park \cite{par2} expressed the joint entropy of $\displaystyle X_{(1:n)},\cdots,X_{(r:n)}$, $H_{1\cdots r:n}$, in the form $\displaystyle H_{1\cdots r:n}=-\ln\frac{n!}{(n-r)!}+n\bar{H}_{1\cdots r:n},$ where
\begin{equation*}\label{hbar1r}
\bar{H}_{1\cdots r:n}=-E\left(\int_0^{U_{(r:n-1)}}\ln(\frac{\mathrm{d}F^{-1}(p)}{\mathrm{d}p})\mathrm{d}p\right)-E\left(\left(1-U_{(r:n-1)}\right)\ln(1-U_{(r:n-1)})\right),
\end{equation*}
and provided its estimator as:
\begin{equation}\label{akram}
\bar{H}_{m,n,r}=\frac{1}{n}\sum_{i=1}^r\ln\left(\frac{x_{(i+m:n)}-x_{(i-m:n)}}{\frac{2m}{n}}\right)-\left(1-\frac{r}{n}\right)\ln\left(1-\frac{r}{n}\right).
\end{equation}
\par
By approximating $\displaystyle \bar{H}_{1\cdots r:n}$ with
\begin{equation}\label{dec}
-\int_0^{\frac{r}{n}}\ln(\frac{\mathrm{d}F^{-1}(p)}{\mathrm{d}p})\mathrm{d}p-(1-\frac{r}{n})\ln(1-\frac{r}{n}),
\end{equation}
we obtain an estimator for (\ref{dec}) as:
\begin{equation}\label{2test}
\bar{H}^{(1)}_{m,n,r}=\frac{1}{n}\sum_{i=1}^{r}\ln \left(\frac{\bar{x}^h_{i-1+m}-\bar{x}^h_{i-1-m}}{\frac{m}{n}}\right)-(1-\frac{r}{n})\ln(1-\frac{r}{n}),
\end{equation}
where $\bar{x}^h_{i-1+m}$ is the harmonic mean of $x_{(i:n)},\cdots,x_{(i-1+m:n)}$ and $\bar{x}^h_{i-1-m}$ is the harmonic mean of $x_{(i-1-m:n)},\cdots,x_{(i:n)}$ and the window size $m$ is a positive integer, which is less than $r/2$; and $x_{(i:n)}=x_{(1:n)}$ for $i<1$, and $x_{(i:n)}=x_{(r:n)}$ for $i<r.$ We expect that the performance of this estimator is better than (\ref{akram}), because we use more information for its calculation.

\par
We can easily prove that the scale of the random variable $X$ has no effect on the accuracy of $\bar{H}^{(1)}_{m,n,r}$ in estimating $H_{1\cdots r:n}$.
\begin{property}
Let $H_{1\cdots r:n}^Y$ and $H_{1\cdots r:n}^W$ denote entropies of the distribution of continuous random variables $Y$ and $W$, respectively, and $W=kY$, where $k>0$. It is easy to see that $\displaystyle\bar{x}^{h,W}_{j}=k\bar{x}^{h,Y}_{j}$ for $i=1,\cdots, r$. So we have $\displaystyle\bar{H}^{{(1)}\;W}_{m,n,r}=\frac{r}{n}\ln k+\bar{H}^{{(1)}\;Y}_{m,n,r}.$ Then the following properties hold
\begin{itemize}
\item {$E(\bar{H}^{{(1)}\;W}_{m,n,r})= E(\bar{H}^{{(1)}\;Y}_{m,n,r})+\frac{r}{n}\ln k,$}
\item {$Var(\bar{H}^{{(1)}\;W}_{m,n,r})= Var(\bar{H}^{{(1)}\;Y}_{m,n,r}),$}
\item {$MSE(\bar{H}^{{(1)}\;W}_{m,n,r})= MSE(\bar{H}^{{(1)}\;Y}_{m,n,r}),$}
\end{itemize}
where the superscript $Y$ and $W$ refer to the corresponding distribution.
\end{property}

\begin{lemma}
If $m, n\rightarrow\infty$ and $\frac{m}{n}\rightarrow0$, then $\bar{H}^{(1)}_{m,n,r}-\bar{H}_{m,n,r}\rightarrow0$, which $\bar{H}_{m,n,r}$ is defined in (\ref{akram}).
\end{lemma}
{\bf\em {Proof:}} If we prove that $\left|\bar{H}^{(1)}_{m,n,r}-\bar{H}_{m,n,r}\right|\rightarrow0$ then by Squeeze theorem, $\bar{H}^{(1)}_{m,n,r}-\bar{H}_{m,n,r}\rightarrow0$. So we establish $\left|\bar{H}^{(1)}_{m,n,r}-\bar{H}_{m,n,r}\right|\rightarrow0$ as follows:
\begin{eqnarray*}
0\leq\left|\bar{H}^{(1)}_{m,n,r}-\bar{H}_{m,n,r}\right|&=&\left|\frac{1}{n}\sum_{i=1}^{r}\ln\frac{\bar{x}^h_{i-1+m}-\bar{x}^h_{i-1-m}}{2(x_{(i+m:n)}-x_{(i-m:n)})}\right|\\
&\leq&\frac{1}{n}\sum_{i=1}^{r}\left|\ln\frac{\bar{x}^h_{i-1+m}-\bar{x}^h_{i-1-m}}{2(x_{(i+m:n)}-x_{(i-m:n)})}\right|\\
&\leq&\frac{1}{n}\sum_{i=1}^{r}\left|\ln\frac{x_{(i-1+m:n)}-x_{(i-1-m:n)}}{2(x_{(i+m:n)}-x_{(i-m:n)})}\right|\\
 &\rightarrow&0, ~\mbox{as}~ m,n\rightarrow0~\mbox{and}~\frac{m}{n}\rightarrow0.
\end{eqnarray*}
The first inequality arises by using the Triangle inequality, and the second inequality is true because
\begin{equation*}
\bar{x}^h_{i-1+m}\leq x_{(i-1+m:n)}~\mbox{and}~\bar{x}^h_{i-1+m}\geq x_{(i:n)},
\end{equation*}
therefore, $\left|\bar{H}^{(1)}_{m,n,r}-\bar{H}_{m,n,r}\right|\rightarrow0$. This completes the proof. $\blacksquare$

\begin{theorem}
If $m, n\rightarrow\infty$ and $\frac{m}{n}\rightarrow0$, then $\bar{H}^{(1)}_{m,n,r}$ is a consistent estimator of $\bar{H}_{1\cdots r:n}$.
\end{theorem}
{\bf\em {Proof:}} Park \cite{par2} showed that $\bar{H}_{m,n,r}$ is a consistent estimator of $\bar{H}_{1\cdots r:n}$. So
\begin{eqnarray}
\label{1a} E\left[\bar{H}_{m,n,r}\right] &\rightarrow&  \bar{H}_{1\cdots r:n},\\
\label{4a} Var\left[\bar{H}_{m,n,r}\right] &\rightarrow& 0,
\end{eqnarray}
as $m, n\rightarrow\infty$, and $\frac{m}{n}\rightarrow 0$.
According to the previous Lemma, $\bar{H}^{(1)}_{m,n,r}-\bar{H}_{m,n,r}\rightarrow 0$, so (Billingsley \cite{bil})
\begin{eqnarray}
\label{2a}    E\left[\bar{H}^{(1)}_{m,n,r}-\bar{H}_{m,n,r}\right] &\rightarrow&  0,\\
\label{3a}    E\left[\bar{H}^{(1)}_{m,n,r}-\bar{H}_{m,n,r}\right]^2 &\rightarrow&  0.
\end{eqnarray}
Now, by using (\ref{1a}) and (\ref{2a}), we conclude $E[\bar{H}^{(1)}_{m,n,r}]\rightarrow \bar{H}_{1\cdots r:n}$ (Billingsley \cite{bil}).\\
On the other hand, using (\ref{2a}) and (\ref{3a}), we have
\begin{equation}
\label{5a}
Var[\bar{H}^{(1)}_{m,n,r}-\bar{H}_{m,n,r}]=E[\bar{H}^{(1)}_{m,n,r}-\bar{H}_{m,n,r}]^2-E^2[\bar{H}^{(1)}_{m,n,r}-\bar{H}_{m,n,r}]\rightarrow0,
\end{equation}
Also,
\begin{equation}\label{end}
Var[\bar{H}^{(1)}_{m,n,r}-\bar{H}_{m,n,r}]=Var[\bar{H}^{(1)}_{m,n,r}]+Var[\bar{H}_{m,n,r}]-2Cov(\bar{H}^{(1)}_{m,n,r},\bar{H}_{m,n,r}),
\end{equation}
and by using (\ref{4a})
\begin{equation}
\label{6a}Var[\bar{H}_{m,n,r}]-2Cov(\bar{H}^{(1)}_{m,n,r},\bar{H}_{m,n,r})\rightarrow0.
\end{equation}
So, by applying (\ref{5a}), (\ref{end}) and (\ref{6a}), we deduce $Var[\bar{H}^{(1)}_{m,n,r}]\rightarrow 0$. Therefore,
\begin{eqnarray*}
\label{2}    E\left[\bar{H}^{(1)}_{m,n,r}\right] &\rightarrow&  \bar{H}_{1\cdots r:n}, \\
\label{3}    Var\left[\bar{H}^{(1)}_{m,n,r}\right] &\rightarrow&  0, ~~\mbox{as}~~m,n\rightarrow\infty,~\frac{m}{n}\rightarrow0,
\end{eqnarray*}
so $\bar{H}^{(1)}_{m,n,r}$ is a consistent estimator of $\bar{H}_{1\cdots r:n}~\blacksquare.$

\subsection{Entropy estimator for monotone increasing hazard function alternatives}
\subsubsection{Moving average method}
\par
In statistics, smoothing a data set is   to create an approximating function that attempts to capture important patterns in the data, while leaving out noise phenomena. One of the most common smoothing methods is moving average. This method is a technique that can be applied for the time series analysis, either to produce smoothed periodogram of data, or to make better estimation and forecasts \cite{bro}.
\par
A moving average (MA) method is the unweighted mean of the previous $n$ datum points. Suppose individual observations, $X_1,\cdots,X_n$ are collected. The moving average of width $w$ at time $i$ is defined by Montgomery \cite{mon}
$$Y_i=\frac{X_i+X_{i-1}+\cdots+X_{i-w+1}}{w}=\frac{\sum_{j=i-w+1}^iX_{j}}{w}\quad i\geq w.$$
For periods $i<w$, we do not have $w$ observations to calculate a moving average of width $w$.
\par
Now, we develop the construction of the moving average method. For this aim, we defined the moving average of width $w$ at time $i$ as:
\begin{equation}\label{rev}
Y_i=\frac{X_i+X_{i+1}+\cdots+X_{i+w-1}}{w}$$$$\hspace{.4in}=\frac{\sum_{j=i}^{i+w-1}X_{j}}{w}\quad i\leq n-w+1.
\end{equation}
From Equation (\ref{rev}), the moving average statistic is the average of the $w$ most recent observations. However, for $i > n-w+1$, the moving average at time $i$ is defined as the average of all observations that are equal or greater than $X_{i}$, i.e.
\begin{equation}\label{rev2}
Y_i=\frac{\sum_{j=i}^{n}X_{j}}{n-i+1}\quad i> n-w+1.
\end{equation}
\par
One characteristic of the MA is that if the data have an uneven path, applying the MA will eliminate abrupt variation and cause the smooth path. In the next subsection, this characteristic of the MA method is used and a new entropy estimator is presented.

\subsubsection{Entropy estimator}
\par
In this subsection, we use the MA method and obtain an entropy estimator which provides a new test statistic that achieves higher power than the previous test statistics against the alternatives with monotone increasing hazard functions. According to the subsection 2.1, we know that
\begin{equation*}
\bar{H}_{1\cdots r:n}=-E\left(\int_0^{U_{(r:n-1)}}\ln(\frac{\mathrm{d}F^{-1}(p)}{\mathrm{d}p})\mathrm{d}p\right)-E\left(\left(1-U_{(r:n-1)}\right)\ln(1-U_{(r:n-1)})\right),
\end{equation*}
and the approximated of it, is defined in (\ref{dec}).
\par
$F^{-1}(p)$ as a function of quantiles in (\ref{dec}) is the sample path of order statistics, but usually it is not smooth. So we propose to imply the MA method of proper order, say $k$, to smooth this sample path and define the new variables $y_1,\cdots,y_r$ from the equation (\ref{rev}) and (\ref{rev2}) as follows:
\begin{eqnarray}
\label{y}
 y_{1}&=&\frac{x_{(1:r)}+\cdots +x_{(k:r)}}{k}, \nonumber\\
 y_{2}&=&\frac{x_{(2:r)}+\cdots +x_{(k+1:r)}}{k}, \nonumber\\
\vdots \\
y_{r-k+1}&=&\frac{x_{(r-k+1:r)}+\cdots +x_{(r:r)}}{k}, \nonumber\\
  y_{r-k+2}&=&\frac{x_{(r-k+2:r)}+\cdots +x_{(r:r)}}{k-1},\nonumber \\
 \vdots \nonumber\\
 y_{r-1}&=&\frac{x_{(r-1:r)}+x_{(r:r)}}{2}, \nonumber\\
  y_{r}&=&x_{(r:r)}\nonumber.
\end{eqnarray}

By this method, we obtain an estimator for (\ref{dec}) as:
\begin{equation}
\label{hbar}
\bar{H}^{(2)}_{m,n,r}=\frac{1}{n}\sum_{i=1}^r\ln\frac{y_{(i+m:n)}-y_{(i-m:n)}}{\hat{F}_n(y_{(i+m:n)})-\hat{F}_n(y_{(i-m:n)})}
-\left(1-\frac{r}{n}\right)\ln\left(1-\frac{r}{n}\right),
\end{equation}
where the window size of $m$ is a positive integer, which is less than $r/2+k$; and $x_{(i:n)}=x_{(1:n)}$ for $i<1$, and $x_{(i:n)}=x_{(r:n)}$ for $i<r.$ Also $\hat{F}_n(y_{(i:n)})$ was introduced by Yousefzadeh and Arghami \cite{you} as:
$$\hat{F}_n(y_{(i:n)})=\frac{r-1}{r(n+1)}\left(i+\frac{1}{r-1}+\frac{y_{(i:n)}-y_{(i-1:n)}}{y_{(i+1:n)}-y_{(i-1:n)}}\right),\;\;\; i=1,\cdots,r,$$
for $y<y_{(1:n)},\;\hat{F}_n(y)$ is less than $\frac{1}{n+1}$ and for $y>y_{(r:n)},\;\hat{F}_n(y)$ is more than $\frac{r}{n+1}$.
\par
We can prove that the scale of the random variable $X$ has no effect on the accuracy of $\bar{H}^{(2)}_{m,n,r}$ in estimating $H_{1\cdots r:n}$.

\begin{property}
Let $H_{1\cdots r:n}^Y$ and $H_{1\cdots r:n}^W$ denote entropies of the distribution of continuous random variables $Y$ and $W$, respectively, and $W=kY$, where $k>0$. It is easy to see that
$$\hat{F}^Y_n(w_{(i:n)})=\frac{r-1}{r(n+1)}\left(i+\frac{1}{r-1}+\frac{w_{(i:n)}-w_{(i-1:n)}}{w_{(i+1:n)}-w_{(i-1:n)}}\right)=\hat{F}^W_n(y_{(i:n)})$$
for $i=1,\cdots, r$. So we have
\begin{eqnarray*}
\bar{H}^{(2)\;W}_{m,n,r}&=&\frac{1}{n}\sum_{i=1}^r\ln\frac{ky_{(i+m:n)}-ky_{(i-m:n)}}{\hat{F}_n(ky_{(i+m:n)})-\hat{F}_n(ky_{(i-m:n)})}
\\ &-& \left(1-\frac{r}{n}\right)\ln\left(1-\frac{r}{n}\right)=\frac{r}{n}\ln k+\bar{H}^{(2)\;Y}_{m,n,r}.
\end{eqnarray*}
Then the following properties hold
\begin{itemize}
\item {$E(\bar{H}^{{(2)}\;W}_{m,n,r})= E(\bar{H}^{{(2)}\;Y}_{m,n,r})+\frac{r}{n}\ln k,$}
\item {$Var(\bar{H}^{{(2)}\;W}_{m,n,r})= Var(\bar{H}^{{(2)}\;Y}_{m,n,r}),$}
\item {$MSE(\bar{H}^{{(2)}\;W}_{m,n,r})= MSE(\bar{H}^{{(2)}\;Y}_{m,n,r}),$}
\end{itemize}
where the superscript $Y$ and $W$ refer to the corresponding distribution.
\end{property}

\begin{example}
For the explanation of the proposed method, we simulate 30 samples from the exponential distribution with mean $1$, consider their order statistics and censor 5 of them from the right, and plot the sample path of 25 points in Figure $1$ with $I$.
\input{epsf}
\epsfxsize=5in \epsfysize=4in
\begin{figure}
\label{order}
\centerline{\epsfxsize=4in \epsfysize=3in \epsffile{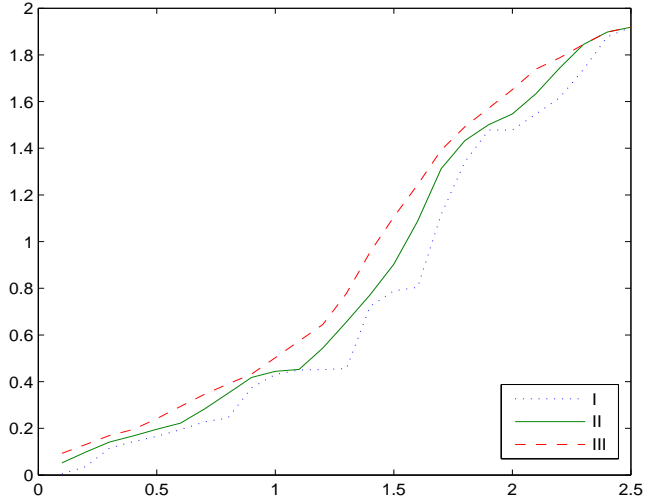}}
\vspace{-0.3in}
\caption{\footnotesize Sample path of order statistics from the exponential distribution (I), smoothed path of order 3 (II) and order 5 (III) at MA method. }
\end{figure}
The sample path of order statistics is smoothed by MA of order $3$. New variables are defined from (\ref{y}) and the smoothed path of new variables is plotted in Figure 1 with $II$. This plot shows that the new sample path is smoother than the sample path of the original order statistics. Also, with considering MA of order 5, we define new variables from (\ref{y}) and plot them in Figure 1 with $III$. Even though the smoothing sample path of order statistics by using MA of order 3 is not as smooth as using MA of order 5, the resulting powers, which are demonstrated in section $4$, are the same up to two digits of decimals. So without loss of generality, we just consider MA of order $k=3$ in (\ref{y}).
\end{example}

\begin{lemma}
If $m, n\rightarrow\infty$ and $\frac{m}{n}\rightarrow0$, then $\bar{H}^{(2)}_{m,n,r}-\bar{H}_{m,n,r}\rightarrow0$.
\end{lemma}
{\bf\em {Proof:}} If we prove that $\left|\bar{H}^{(2)}_{m,n,r}-\bar{H}_{m,n,r}\right|\rightarrow0$ then by Squeeze theorem, $\bar{H}^{(2)}_{m,n,r}-\bar{H}_{m,n,r}\rightarrow0$. So we establish $\left|\bar{H}^{(2)}_{m,n,r}-\bar{H}_{m,n,r}\right|\rightarrow0$ as follows:
\begin{eqnarray*}
0\leq\left|\bar{H}^{(2)}_{m,n,r}-\bar{H}_{m,n,r}\right|&=&\left|\frac{1}{n}\sum_{i=1}^{r}\ln\frac{\frac{2m}{n}(y_{(i+m:n)}-y_{(i-m:n)})}{\left(\hat{F_n}(y_{(i+m:n)})-\hat{F_n}
(y_{(i-m:n)})\right)(x_{(i+m:n)}-x_{(i-m:n)})}\right|\\
&\leq&\frac{1}{n}\sum_{i=1}^{r}\left|\ln\frac{\frac{2m}{n}(y_{(i+m:n)}-y_{(i-m:n)})}{\left(\hat{F_n}(y_{(i+m:n)})-\hat{F_n}
(y_{(i-m:n)})\right)(x_{(i+m:n)}-x_{(i-m:n)})}\right|\\
&\leq&\frac{1}{n}\sum_{i=1}^{r}\left|\ln\left[\frac{n+1}{2m-1}\frac{2m}{n}\frac{y_{(i+m:n)}-y_{(i-m:n)}}{x_{(i+m:n)}-x_{(i-m:n)}}\right]\right|\\
&\leq&\frac{1}{n}\sum_{i=1}^{r}\left|\ln\left[\frac{n+1}{2m-1}\frac{2m}{n}\frac{x_{(i+m+k:n)}-x_{(i-m:n)}}{x_{(i+m:n)}-x_{(i-m:n)}}\right]\right|\\
 &\rightarrow&0, ~\mbox{as}~ m,n\rightarrow0~\mbox{and}~\frac{m}{n}\rightarrow0.
\end{eqnarray*}
The first inequality arises by using the Triangle inequality and the second inequality is true, for more details see Yousefzadeh and Arghami \cite{you}. Also, the third inequality is true because
\begin{eqnarray*}
y_{(i+m:n)}=\frac{1}{k}\left(x_{(i+m:n)}+\cdots+x_{(i+m+k:n)}\right)&\Rightarrow& y_{(i+m:n)}\leq x_{(i+m+k:n)}\\
y_{(i-m:n)}=\frac{1}{k}\left(x_{(i-m:n)}+\cdots+x_{(i-m+k:n)}\right)&\Rightarrow& y_{(i-m:n)}\geq x_{(i-m:n)},
\end{eqnarray*}
therefore, $\left|\bar{H}^{(2)}_{m,n,r}-\bar{H}_{m,n,r}\right|\rightarrow0$. This completes the proof. $\blacksquare$

\begin{theorem}
If $m, n\rightarrow\infty$ and $\frac{m}{n}\rightarrow0$, then $\bar{H}^{(2)}_{m,n,r}$ is a consistent estimator of $\bar{H}_{1\cdots r:n}$.
\end{theorem}
{\bf\em {Proof:}} We should prove that
\begin{eqnarray*}
\label{2}    E\left[\bar{H}^{(2)}_{m,n,r}\right] &\rightarrow&  \bar{H}_{1\cdots r:n}, \\
\label{3}    Var\left[\bar{H}^{(2)}_{m,n,r}\right] &\rightarrow&  0, ~~\mbox{as}~~m,n\rightarrow\infty,~\frac{m}{n}\rightarrow0
\end{eqnarray*}
These equation obtain from the consistency of $\bar{H}_{m,n,r}$ for $\bar{H}_{1\cdots r:n}$. Proof of this theorem is quite similar to the proof of Theorem 2.1. $\blacksquare$

\section{Test statistics}
For a null distribution function $f^0(z;\theta),$ the KL information for the Type-II censored data is defined to be:
$$I_{1\cdots r:n}(f,f^0)=\int_{-\infty}^{\infty}f_{1\cdots r:n}(z;\theta)\ln\frac{f_{1\cdots r:n}(z;\theta)}{f^0_{1\cdots r:n}(z;\theta)}\mathrm{d}z.$$
Then the KL information can be approximated with
\begin{equation}\label{I}
I_{1\cdots r:n}(f,f^0)=-n\bar{H}_{1\cdots r:n}-\sum^{r}_{i=1}\ln f^0(z_{(i:n)};\theta)-(n-r)\ln(1-F^0(z_{(r:n)};\theta)).
\end{equation}
Thus the test statistic based on $I_{1\cdots r:n}(f,f^0)/n$ can be written as:
\begin{equation*}
\label{rtestt}
T^{(j)}_{m,n,r}=-\bar{H}^{(j)}_{m,n,r}-\frac{1}{n}\left(\sum^{r}_{i=1}\ln f^0(z_{(i:n)}^{(j)};\hat{\theta})+(n-r)\ln(1-F^0(z_{(r:n)}^{(j)};\hat{\theta}))\right),~~~~j=1,2,
\end{equation*}
where
\begin{equation}
\label{mat}
z_{(i:n)}^{(j)}=\left\{
\begin{array}{rl}
x_{(i:n)} & {j=1}\\
y_{(i:n)} & {j=2}
\end{array}
\right.
,~~~~i=1\cdots r
\end{equation}
and $\hat{\theta}$ is an estimator of $\theta$ and $\bar{H}^{(1)}_{m,n,r}$ and $\bar{H}^{(2)}_{m,n,r}$ is defined in (\ref{2test}) and (\ref{hbar}), respectively.

\section{Testing exponentiality based on the Kullback-Leibler information}
\subsection{Test statistics}
\par
Suppose that we are interested in a goodness of fit test for
\begin{equation*}
\left\{
\begin{array}{c}
H_0: f^0(x)=\frac{1}{\theta}\exp(-\frac{x}{\theta}),\\
H_1: f^0(x)\neq\frac{1}{\theta}\exp(-\frac{x}{\theta}),
\end{array}
\right.
\end{equation*}
where $\theta$ is unknown. Then the KL information for the Type-II censored data can be approximated in view of (\ref{I}) with
\begin{equation*}
I_{1\cdots r:n}(f;f^0)=-n\bar{H}_{1\cdots r:n}+r\ln\theta+\frac{1}{\theta}\left(\sum_{i=1}^rZ_{(i:n)}^{(j)}+(n-r)Z_{r:n}^{(j)}\right),
\end{equation*}
If we estimate the unknown $\theta$ with the maximum likelihood estimator, $$\displaystyle\hat{\theta}=\left(\sum_{i=1}^rZ_{(i:n)}^{(j)}+(n-r)Z_{(r:n)}^{(j)}\right)/r,$$ then we have two estimators of $I_{1\cdots r:n}(f;f^0)/n$ as:
\begin{equation*}\label{T1}
T^{(j)}_{m,n,r}=-\bar{H}^{(j)}_{m,n,r}+\frac{r}{n}\left\{\ln\left[\frac{1}{r}\left(\sum_{i=1}^rZ_{(i:n)}^{(j)}+(n-r)Z_{(r:n)}^{(j)}\right)\right]+1\right\},
~~~~j=1,2,
\end{equation*}
where the random variable $Z_{(i:n)}^{(j)}$ takes the value $z_{(i:n)}^{(j)}$ which is defined in (\ref{mat}). Since $I$ is non-negative and is zero if and only if $f=f^0$, a.e., we reject the null hypothesis for large values of $T^{(1)}_{m,n,r}$ and $T^{(2)}_{m,n,r}$.

\subsection{Implementation of the test}
\par
Because the sampling distributions of the test statistics are intractable, we determine the percentage points using $10000$ Monte Carlo samples from an exponential distribution. In determining the window size $m$ which depends on $n,\;r$ and the $\alpha$, we define the optimal window size $m$ to be one which gives minimum critical points in the sense of Ebrahimi et al. \cite{ebr2}. However, we find from the simulated percentage points, the optimal window size $m$. In view of these results, our recommended values of $m$ for different $r$ and test statistic $T^{(1)}_{m,n,r}$ are listed in Table 1 and the critical values of $T^{(1)}_{m,n,r}$ corresponding to the optimum values of $m$, are given in Table 2. Also, our recommended values of $m$ for different $r$ and test statistic $T^{(2)}_{m,n,r}$ are listed in Table 3, where $m^*=r/2+3$ and the critical values of $T^{(2)}_{m,n,r}$ corresponding to the optimum values of $m$, are given in Table 4.
\begin{table}
\label{window}
\caption{{\footnotesize Values of the window size $m$ which gives minimum critical values of $\alpha$ less than $0.1$ for $T^{(1)}_{m,n,r}$}}
\begin{center}
{\begin{tabular}{c||ccc}
\hline
\hline
r & 5-19 & 20-40 & 41-50 \\
\hline
m & 3 & 4 & 5 \\
\hline
\end{tabular}}
\end{center}
\end{table}

\begin{table}
\label{cv}
\caption{{\footnotesize Monte carlo estimate of the critical values of $T^{(1)}_{m,n,r}$ where $m$ is determined from Table $1$ }}
\begin{center}
{\begin{tabular}{c|c|c|c|c}
\hline
\hline
n & r & $\alpha=0.1$ & $\alpha=0.05$ & $\alpha=0.025$\\
\hline
10 & 5  & 0.5962& 0.6855& 0.7692 \\
   & 6  & 0.6155& 0.7185& 0.8039 \\
   & 7  & 0.6398& 0.7333& 0.8185 \\
   & 8  & 0.6676& 0.7607& 0.8648 \\
   & 9  & 0.7152& 0.8075& 0.9025  \\
\hline
20 & 10 & 0.3148& 0.3640& 0.4061  \\
   & 11 & 0.3188& 0.3689& 0.4148  \\
   & 12 & 0.3285& 0.3727& 0.4183  \\
   & 13 & 0.3374& 0.3825& 0.4329  \\
   & 14 & 0.3442& 0.3911& 0.4371  \\
   & 15 & 0.3613& 0.4113& 0.4587  \\
   & 16 & 0.3677& 0.4157& 0.4634  \\
   & 17 & 0.3830& 0.4333& 0.4795  \\
   & 18 & 0.4022& 0.4521& 0.5024  \\
   & 19 & 0.4223& 0.4717& 0.5172  \\
\hline
30 & 15 & 0.2239& 0.2599& 0.2904  \\
   & 16 & 0.2293& 0.2625& 0.2922  \\
   & 17 & 0.2320& 0.2659& 0.2945  \\
   & 18 & 0.2376& 0.2707& 0.3009  \\
   & 19 & 0.2425& 0.2757& 0.3077  \\
   & 20 & 0.2470& 0.2800& 0.3108  \\
   & 21 & 0.2537& 0.2834& 0.3121  \\
   & 22 & 0.2568& 0.2918& 0.3259  \\
   & 23 & 0.2634& 0.2949& 0.3272  \\
   & 24 & 0.2691& 0.3017& 0.3318  \\
   & 25 & 0.2736& 0.3090& 0.3398  \\
   & 26 & 0.2822& 0.3136& 0.3488  \\
   & 27 & 0.2910& 0.3239& 0.3538  \\
\hline
\end{tabular}}
\end{center}
\end{table}

\begin{table}
\label{window}
\caption{{\footnotesize Values of the window size $m$ which gives minimum critical values of $\alpha$ less than $0.1$ for $T^{(2)}_{m,n,r}$}}
\begin{center}
{\begin{tabular}{c||cccccc}
\hline
\hline
r & 4-5 & 6-7 & 8-9 & 10-11 & 12-13 & r-(r+1) \\
\hline
m & 5 & 6 & 7 & 8 & 9 & $m^*$ \mbox{(for even $r$)} \\
\hline
\end{tabular}}
\end{center}
\end{table}

\begin{table}
\label{cv}
\caption{{\footnotesize Monte carlo estimate of the critical values of $T^{(2)}_{m,n,r}$ where $m$ is determined from Table $3$ }}
\begin{center}
{\begin{tabular}{c|c|c|c|c}
\hline
\hline
n & r & $\alpha=0.1$ & $\alpha=0.05$ & $\alpha=0.025$\\
\hline
10 & 5  & 0.3445 & 0.4253 & 0.5087 \\
   & 6  & 0.3251 & 0.4128 & 0.5026 \\
   & 7  & 0.3136 & 0.4099 & 0.4929 \\
   & 8  & 0.3104 & 0.4046 & 0.4915 \\
   & 9  & 0.3101 & 0.4038 & 0.4902 \\
\hline
20 & 10 & 0.1474 & 0.1913 & 0.2310 \\
   & 11 & 0.1426 & 0.1830 & 0.2229 \\
   & 12 & 0.1408 & 0.1828 & 0.2197 \\
   & 13 & 0.1369 & 0.1805 & 0.2181 \\
   & 14 & 0.1322 & 0.1773 & 0.2168 \\
   & 15 & 0.1294 & 0.1740 & 0.2160 \\
   & 16 & 0.1281 & 0.1742 & 0.2161 \\
   & 17 & 0.1280 & 0.1739 & 0.2073 \\
   & 18 & 0.1268 & 0.1649 & 0.2072 \\
   & 19 & 0.1200 & 0.1620 & 0.1988 \\
\hline
30 & 15 & 0.0865 & 0.1168 & 0.1500 \\
   & 16 & 0.0859 & 0.1152 & 0.1430 \\
   & 17 & 0.0843 & 0.1124 & 0.1383 \\
   & 18 & 0.0829 & 0.1090 & 0.1380 \\
   & 19 & 0.0824 & 0.1083 & 0.1355 \\
   & 20 & 0.0815 & 0.1079 & 0.1311 \\
   & 21 & 0.0806 & 0.1043 & 0.1294 \\
   & 22 & 0.0777 & 0.1038 & 0.1281 \\
   & 23 & 0.0759 & 0.1027 & 0.1280 \\
   & 24 & 0.0741 & 0.0988 & 0.1275 \\
   & 25 & 0.0699 & 0.0976 & 0.1265 \\
   & 26 & 0.0662 & 0.0977 & 0.1219 \\
   & 27 & 0.0635 & 0.0910 & 0.1200 \\
\hline
\end{tabular}}
\end{center}
\end{table}

\subsection{Power results}
\par
There are lots of test statistics for exponentiality concerning uncensored data including \cite{asc}, \cite{gan}, \cite{hen}-\cite{lar}, but only some of them can be extended to the censored data. We consider here the test statistics of \cite{bra}, and \cite{par2} among them. Brain and Shapiro \cite{bra} proposed two test statistics as:
\begin{eqnarray*}
z&=&\left(\frac{12}{r-2}\right)^\frac{1}{2}\frac{\sum_{i=1}^{r-1}(i-\frac{r}{2})Y_{i+1}}{\sum_{i=1}^{r-1}Y_{i+1}}\\
Z&=&z^2+\left(\frac{5}{4(r+1)(r-2)(r-3)}\right)^\frac{1}{2}\\
&\quad&\times\frac{12\sum_{i=1}^{r-1}(i-\frac{r}{2})^2Y_{i+1}-r(r-2)\sum_{i=1}^{r-1}Y_{i+1}}{\sum_{i=1}^{r-1}Y_{i+1}}
\end{eqnarray*}
where $Y_1=nX_{(1:n)}$, and $Y_i=(n-i+1)(X_{(i:n)}-X_{(i-1:n)})$, $i=2,\cdots,r$; and show that $z$ and $Z$ perform better than other test statistics for the censored data. Recently Park \cite{par2} proposed a test statistic as:
\begin{equation*}
T_{m,n,r}=-\bar{H}_{m,n,r}+\frac{r}{n}\left\{\ln\left[\frac{1}{r}\left(\sum_{i=1}^rX_{(i:n)}+(n-r)X_{(r:n)}\right)\right]+1\right\},
\end{equation*}
where $\bar{H}_{m,n,r}$ is presented in (\ref{akram}). He showed that the power of the proposed test statistic is greater than the power of the test statistics which was introduced by Brain and Shapiro \cite{bra} against the alternatives with monotone increasing hazard functions.
\par
Because the proposed test statistics are essentially related to the hazard function, the alternatives are considered according to the type of hazard functions as follows:
\begin{itemize}
{\item Monotone decreasing hazard: Chi-square with degree of freedom 1 (A1), Gamma with shape parameter 0.5 (A2), Weibull with shape parameter 0.5 (A3) and Generalized Exponential with shape 0.5 (A4).}
{\item Monotone increasing hazard: Uniform (B1), Weibull with shape parameter 2 (B2), Gamma with shape parameter 1.5, 2 (B3, B4 respectively), Chi-square with degree of freedom 3, 4 (B5, B6 respectively), Beta with shape parameters 1 and 2, 2 and 1 (B7, B8 respectively).}
{\item Non-monotone hazard: Log normal with shape parameter 0.6, 1.0, 1.2 (C1, C2, C3 respectively), Beta with shape parameters 0.5 and 1.0 (C4).}
\end{itemize}
We consider here the sample size to be $30$, and draw conclusions. We made 10000 Monte Carlo simulations for $n=30$ to estimate the powers of our proposed test statistics and the competing test statistics, for $\alpha=0.1$. The simulation results are summarized in Figures $2-4$. We can see from these figures that any test statistics does not beat others against all alternatives, but it is notable that the first proposed test statistic, $T^{(1)}_{m,n,r}$, shows better powers than the competing test statistics against the alternatives with monotone decreasing hazard functions, see Figure $2$. Also, against the alternatives with monotone increasing hazard functions, the second proposed test statistic, $T^{(2)}_{m,n,r}$, shows better powers than the competing test statistics, see Figure $3$.

\input{epsf}
\label{f1}
\epsfxsize=5in \epsfysize=4in
\begin{figure}
\vspace{-1in}
\centerline{\epsfxsize=5in \epsfysize=3.5in \epsffile{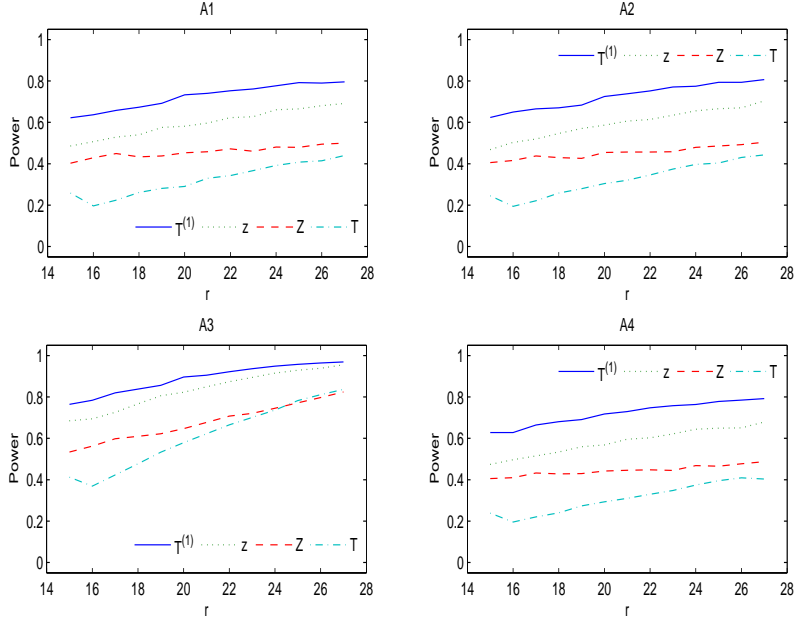}}
\vspace{-0.25in}
\caption{\footnotesize Power comparison: monotone decreasing hazard alternative at 10\% when the sample size is 30. $r$ is the remaining data after the implementation of Type-II censoring scheme. $z$ and $Z$ were introduced by Brain and Shapiro \cite{bra} and $T$ was introduced by Park \cite{par2}. (A1) Chi-square: df 1, (A2) Gamma: shape 0.5, (A3) Weibull: shape 0.5, (A4) Generalized Exponential: shape 0.5.}
\end{figure}

\input{epsf}
\label{f3}
\epsfxsize=5in \epsfysize=4in
\begin{figure}
\vspace{-2in}
\centerline{\epsfxsize=5in \epsfysize=4in \epsffile{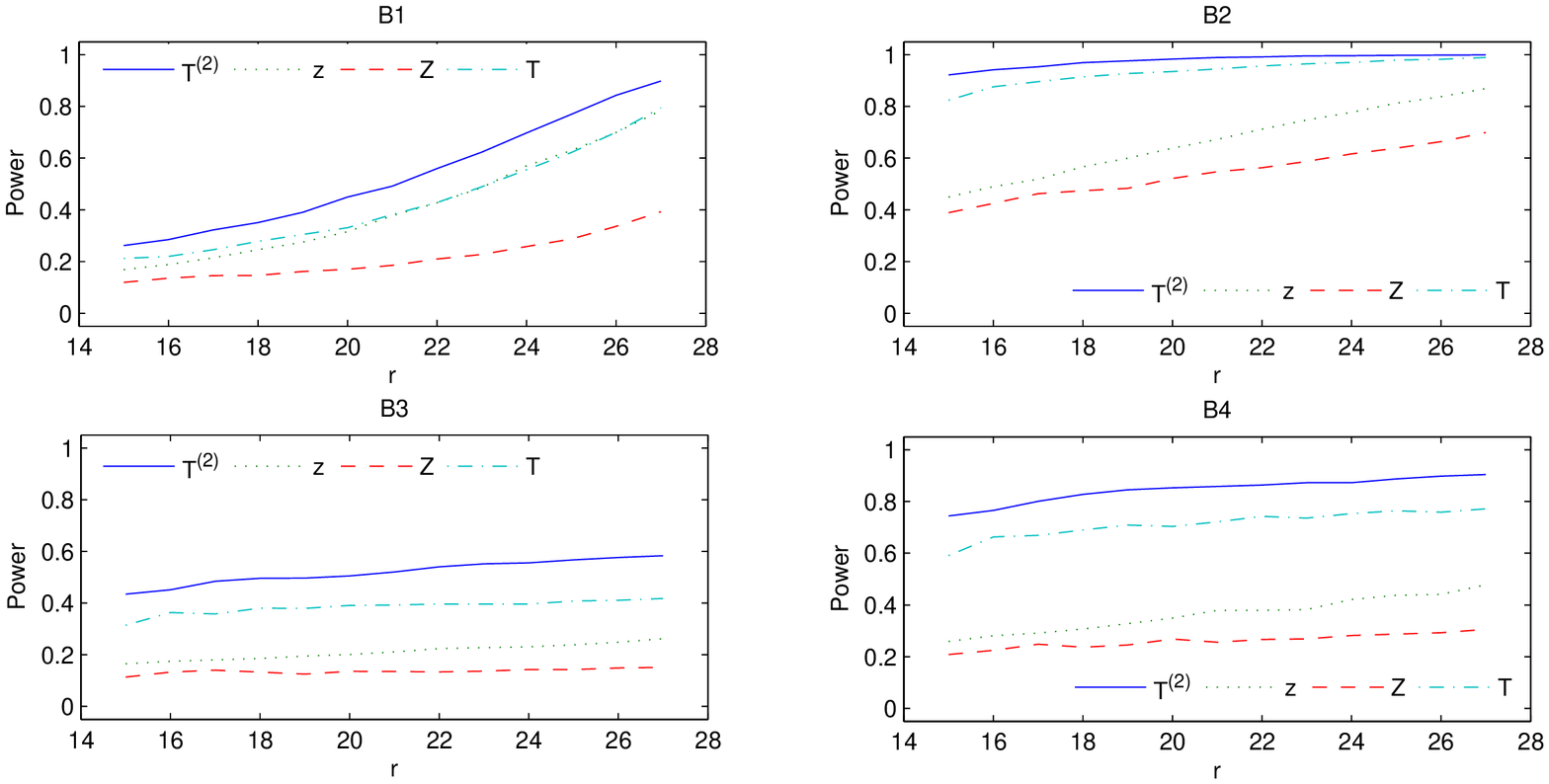}}
\vspace{-1.8in}
\centerline{\epsfxsize=5in \epsfysize=4in \epsffile{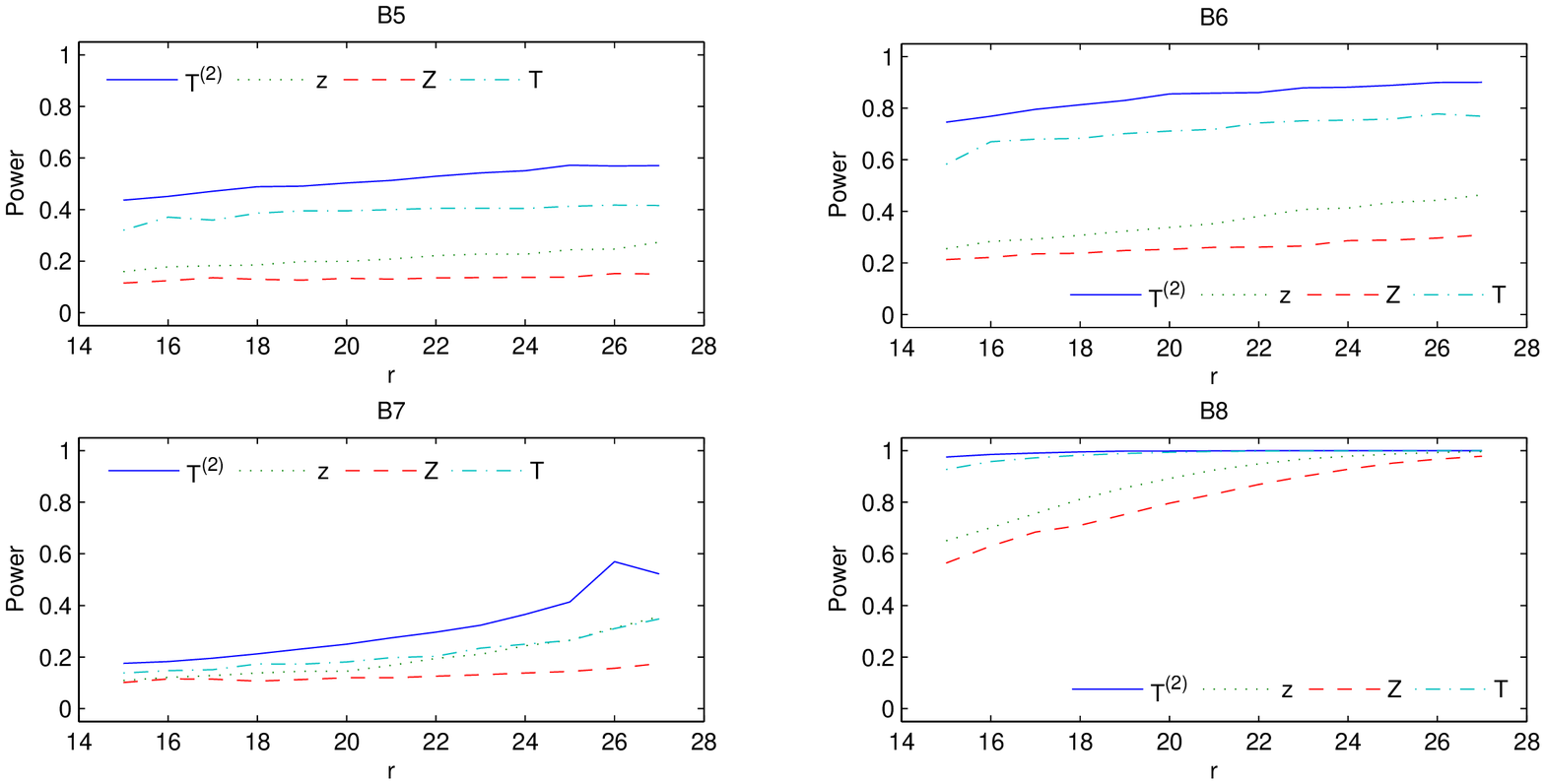}}
\vspace{-.25in}
\caption{\footnotesize Power comparison: monotone increasing hazard alternative at 10\% when the sample size is 30. $r$ is the remaining data after the implementation of Type-II censoring scheme. $z$ and $Z$ were introduced by Brain and Shapiro \cite{bra} and $T$ was introduced by Park \cite{par2}. (B1) Uniform, (B2) Weibull: shape 2, (B3) Gamma: shape 1.5, (B4) Gamma: shape 2, (B5) Chi-square: df 3, (B6) Chi-square: df 4, (B7) Beta: shape 1 and 2, (B8) Beta: shape 2 and 1.}
\end{figure}

\input{epsf}
\label{f2}
\epsfxsize=5in \epsfysize=4in
\begin{figure}
\vspace{-1.2in}
\centerline{\epsfxsize=5in \epsfysize=3.5in \epsffile{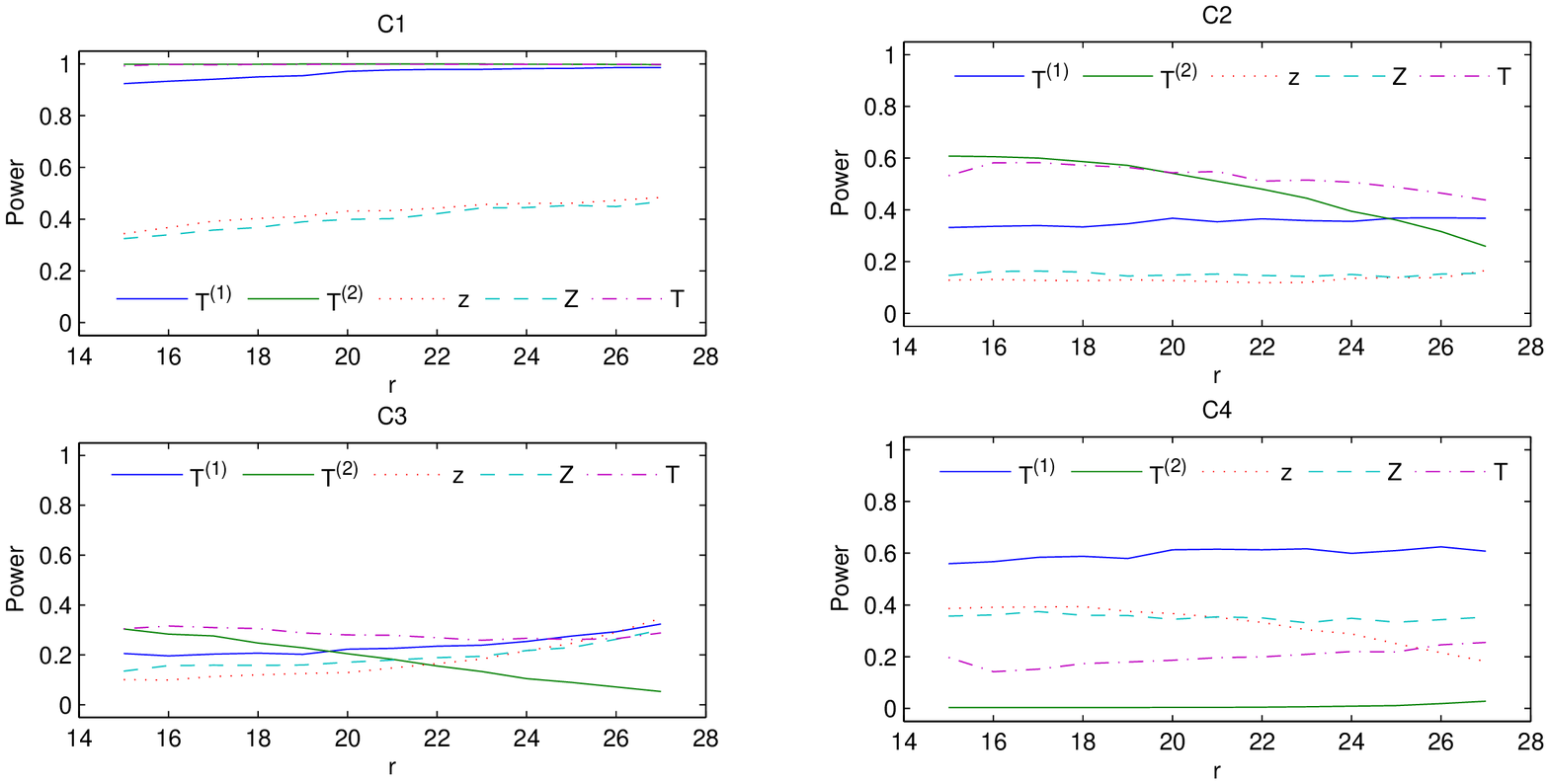}}
\vspace{-0.25in}
\caption{\footnotesize Power comparison: non-monotone hazard alternative at 10\% when the sample size is 30. $r$ is the remaining data after the implementation of Type-II censoring scheme. $z$ and $Z$ were introduced by Brain and Shapiro \cite{bra} and $T$ was introduced by Park \cite{par2}. (C1) Log normal: shape 0.6, (C2) Log normal: shape 1, (C3) Log normal: shape 1.2, (C4) Beta: shape 0.5 and 1.}
\end{figure}

\subsection{RMSE comparisons}
\par
In this subsection, we report the results of a simulation study which compares the performances of the introduced entropy estimators with the estimator proposed by Park \cite{par2} in terms of their biases and RMSEs. We consider here the sample size to be $30$, and draw conclusions. We made 10000 Monte Carlo simulations for $n=30$ and different $r$ to obtain the $\bar{H}_{m,n,r}$, $\bar{H}^{(1)}_{m,n,r}$, $\bar{H}^{(2)}_{m,n,r}$, their biases and RMSEs. The simulation results are summarized in Table $5$. The results show that $\bar{H}^{(2)}_{m,n,r}$ has the smallest bias and RMSE among them. Also, the bias and RMSE of $\bar{H}_{m,n,r}$ is smaller than $\bar{H}^{(1)}_{m,n,r}$. We plot the empirical density of the test statistics based on these estimators for other $n$ and $r$ in Figure $5$. This figure confirms the simulation results.

\begin{table}
\label{cv}
\caption{{\footnotesize Monte Carlo biases and root of mean square errors (RMSE) for exponential distribution}}
\begin{center}
\begin{tabular}{c|c|c|c|c||c|c|c}
\hline
\hline
 \multicolumn{1}{c|}{} & \multicolumn{1}{c|}{} & \multicolumn{3}{c||}{Bias} & \multicolumn{3}{c}{RMSE} \\
\cline{3-8}
n   & r  & $\bar{H}^{(1)}_{m,n,r}$ & $\bar{H}^{(2)}_{m,n,r}$ & $\bar{H}_{m,n,r}$ & $\bar{H}^{(1)}_{m,n,r}$ & $\bar{H}^{(2)}_{m,n,r}$ & $\bar{H}_{m,n,r}$ \\
\hline
 30 & 15 & -0.1626 & -0.0100 & -0.1370 &  0.2159 & 0.1426 & 0.1953 \\
    & 16 & -0.1691 & -0.0102 & -0.1508 &  0.2245 & 0.1478 & 0.2078 \\
    & 17 & -0.1717 & -0.0035 & -0.1521 &  0.2284 & 0.1511 & 0.2108 \\
    & 18 & -0.1760 &  0.0019 & -0.1540 &  0.2348 & 0.1557 & 0.2156 \\
    & 19 & -0.1788 &  0.0095 & -0.1543 &  0.2386 & 0.1594 & 0.2176 \\
    & 20 & -0.1902 &  0.0150 & -0.1579 &  0.2494 & 0.1638 & 0.2233 \\
    & 21 & -0.1957 &  0.0189 & -0.1616 &  0.2565 & 0.1702 & 0.2294 \\
    & 22 & -0.1954 &  0.0316 & -0.1588 &  0.2575 & 0.1752 & 0.2291 \\
    & 23 & -0.2019 &  0.0384 & -0.1630 &  0.2660 & 0.1821 & 0.2359 \\
    & 24 & -0.2042 &  0.0529 & -0.1626 &  0.2683 & 0.1882 & 0.2363 \\
    & 25 & -0.2145 &  0.0613 & -0.1687 &  0.2792 & 0.1960 & 0.2442 \\
    & 26 & -0.2169 &  0.0843 & -0.1665 &  0.2823 & 0.2106 & 0.2440 \\
    & 27 & -0.2300 &  0.0980 & -0.1745 &  0.2941 & 0.2204 & 0.2514 \\
\hline
\end{tabular}
\end{center}
\end{table}

\input{epsf}
\label{f1}
\epsfxsize=5in \epsfysize=4in
\begin{figure}
\centerline{\epsfxsize=6in \epsfysize=3in \epsffile{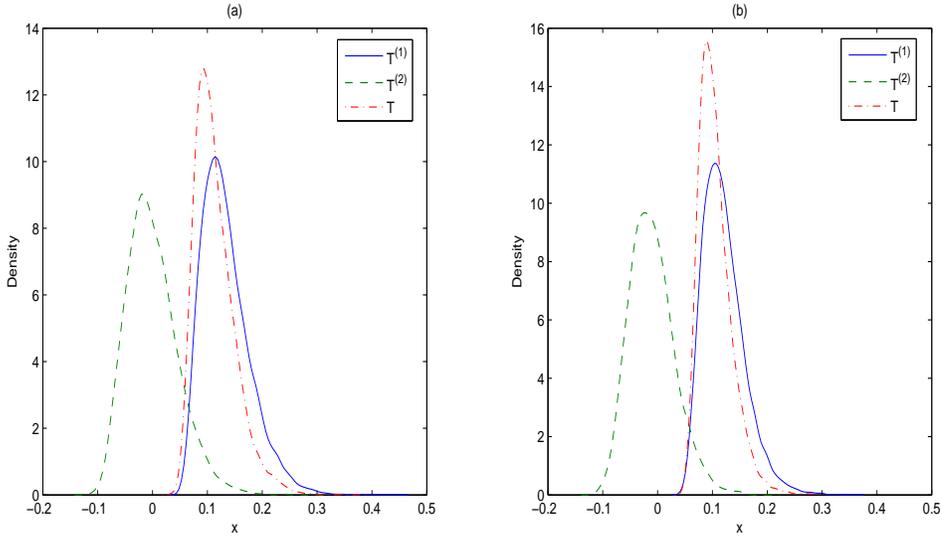}}
\vspace{-0.15in}
\caption{\footnotesize Empirical density functions of $T^{(1)}_{m,n,r}$, $T^{(2)}_{m,n,r}$ and $T_{m,n,r}$ based on 10000 simulations (a) $n=40$ and $r=25$ (b) $n=50$ and $r=35$ under the exponential hypothesis.}
\end{figure}

\section{Conclusion}
\par
In this paper, the entropy estimator of the Type-II censored data which was introduced by Park \cite{par2} is modified and two new entropy estimators are obtained. Simulation results showed that the second proposed entropy estimator compared favourably with their competitors in terms of bias and RMSE, as it is expected of the structure of $\bar{H}^{(2)}_{m,n,t}$. Also, we provided two new test statistics for testing exponentiality with the Type-II censored data. The first one was quite powerful when compared to the existing goodness of fit tests proposed against the alternatives with monotone decreasing hazard functions. Moreover, the second one showed better powers than the available test statistics against the alternatives with monotone increasing hazard functions.
\par
This work has the potential to be applied in the context of censored data and goodness of fit tests. This paper can elaborate further researches by extending such modifications for other censoring schemes such as progressive censoring schemes. Finally, this area of research can be expanded by considering other distributions besides the exponential distribution such as Pareto, Log normal and Weibull distributions.

\bibliographystyle{plain}

\end{document}